\documentclass[journal]{IEEEtran}
\ifCLASSINFOpdf
  \usepackage[pdftex]{graphicx}
\else
\fi
\usepackage[cmex10]{amsmath}
\usepackage{amssymb}  
\usepackage{algorithm}
\usepackage{algorithmic}
\usepackage{subcaption}
\begin{document}
%
\title{Particle Swarm Optimization Based Source Seeking}
%
%
%

\author{Rui Zou, Vijay Kalivarapu, Eliot Winer, James Oliver, Sourabh Bhattacharya%
\thanks{Rui Zou, Eliot Winer, James Oliver and Sourabh Bhattacharya are with the Department of Mechanical Engineering,{\tt\small \{rzou,ewiner,oliver,sbhattac\}@iastate.edu}, Vijay Kalivarapu is with the Virtual Reality Applications Center, {\tt\small vkk2@iastate.edu}, Iowa State University, Ames, IA, 50011 USA.}} 
\maketitle

\begin{abstract}
Signal source seeking using autonomous vehicles is a complex problem. The complexity increases manifold when signal intensities captured by physical sensors onboard are noisy and unreliable. Added to the fact that signal strength decays with distance, noisy environments make it extremely difficult to describe and model a decay function. This paper addresses our work with seeking maximum signal strength in a continuous electromagnetic signal source with mobile robots, using Particle Swarm Optimization (PSO). A one to one correspondence with swarm members in a PSO and physical mobile robots is established and the positions of the robots are iteratively updated as the PSO algorithm proceeds forward. Since physical robots are responsive to swarm position updates, modifications were required to implement the interaction between real robots and the PSO algorithm. The development of modifications necessary to implement PSO on mobile robots, and strategies to adapt to real life environments such as obstacles and collision objects are presented in this paper. Our findings are also validated using experimental testbeds.

\end{abstract}

%

Note to Practitioners:
\begin{abstract}
This paper was inspired by the source seeking problem when the signal source is very noisy and non-smooth. We focus our work particularly on electromagnetic source, but the strategies proposed in this paper can be applied to other types of sources. Most existing strategies approach this problem with gradient-based methods using one or more mobile agents. These methods either assume the signal profiles to be smooth or use complicated and computationally costly procedures to obtain accurate estimation of gradient information. This paper suggests approaching the problem using a heuristic method which is simple to implement on robots with limited computation capability. We propose some modifications to a population based optimization technique – Particle Swarm Optimization (PSO), so as to adapt it a real-world scenario where a group of mobile agents trying to find an unknown electromagnetic source. The mobile agents know their own positions and can measure the signal strength at their current positions. They can share information and plan for the next step based on individual and group memory. We then propose a complete solution to ensure the effectiveness of PSO in complex environments where collisions may occur. We incorporate static and dynamic obstacle avoidance strategies in PSO to make it fully applicable to real-world scenario. In the end, we validate the proposed method in experiments. In our future work, we will work on improving the efficiency of the method.
\end{abstract}


%
\IEEEpeerreviewmaketitle

\section{Introduction}
Seeking a source with autonomous vehicles is an area of growing interest and wide applications. The source could be an electromagnetic signal, acoustic signal, thermal signal, or a chemical/biological agent. Motivated from source-seeking behavior exhibited by natural species from a microscopic level \cite{Optimotaxis} to a macroscopic level \cite{AnimalNavigation}, researchers have developed robots \cite{TDoA} and sensor networks \cite{Coop_control} that can imitate these behaviors in order to perform complex tasks such as environment monitoring, search and rescue operations, explosive detection, drug detection, sensing leakage of hazardous chemicals, pollution sensing and environmental studies.

\begin{figure}[htb]
	\centering
	\begin{subfigure}[b]{\linewidth}
		\includegraphics[width=\textwidth]{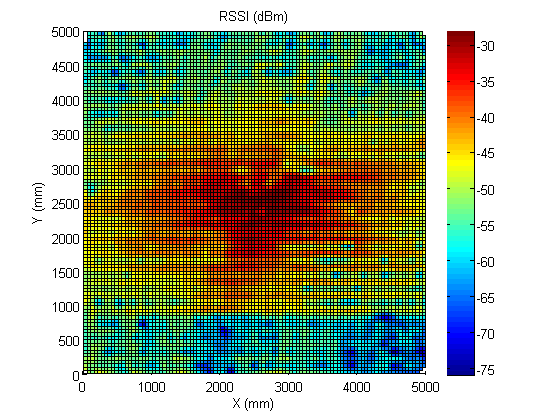}
	\end{subfigure}\\%
	\begin{subfigure}[b]{\linewidth}
		\includegraphics[width=\textwidth]{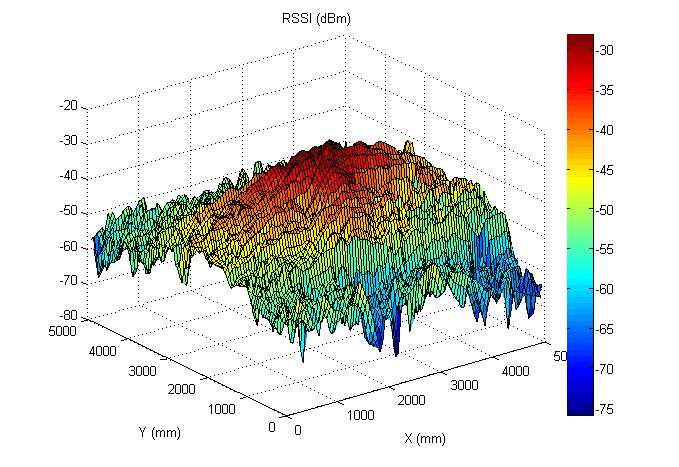}
	\end{subfigure}
	\caption{Map of RSSI (Received Signal Strength Indication). The color bars on the right indicate the signal strength in dBm. In this map, the source is located in the center where the RSSI is approximately -29 dBm}
	\label{fig:Map}
\end{figure}

In this work, we address the problem in which a team of mobile agents, called the seekers, attempt to find the location of a source that emits an electromagnetic signal of unknown strength. The seekers can continuously sense the signal strength transmitted by the source at their current positions which generally decays with distance from the source. The decay profile of the signal strength is very noisy as shown in Figure \ref{fig:Map}, which makes many existing methods inapplicable. Based on this information, we investigate the issue of modifying Particle Swarm Optimization (PSO) and applying it to swarm mobile robots to approach the source seeking problem.
 
A vast amount of research has been done on source seeking with autonomous agents based on the idea of gradient descent/ascent. Many methods employ a single agent to search for single or multiple, static targets. Authors of \cite{Nehorai} propose methods of computing the gradients of the Cramer-Rao bound on the location error with respect to a sensor's coordinates and moving the sensor opposite to the corresponding gradient directions. They extend their techniques in \cite{VaporSource} to provide motion planning strategy for a mobile sensor used for detecting low concentration vapors. This method assumes very small sensing noise in the measurement model, and its applicability to a very noisy model is not known. The same issue is encountered by the algorithm implemented on Autonomous Underwater Vehicles to localize hydrothermal vents in \cite{AUVCircle}. The method of estimating local gradient by taking measurements on a circle may no longer be effective on such a noisy model. In \cite{RotationBasedAngle}, the authors utilize the difference in RSSI received by a rotating directional antenna on the single wing of the Samarai MAV to obtain gradient. However, the application of this method is restricted by the dynamics of the Samarai MAV whose entire body including the mono-wing rotates at all time for stable hovering \cite{SingleWingMAV}. For other types of robots, additional structure has to be made, and extra energy has to be allocated to rotate the antenna at all times which is neither convenient nor cost-efficient.

Recently, extremum seeking \cite{Extremum_seeking} techniques have been adopted aiming at the source seeking problem. It has been applied to nonholonomic vehicles in both 2-D \cite{Non-h_forward_v, Cochran2009} and 3-D \cite{3-D_seeking} environment using sinusoidal and stochastic perturbation \cite{Liu20101443}. The work in \cite{Stanković20101243} also addresses the issue of stochastic noise in measurement. It uses the methodology of stochastic approximation to deal with colored noise. Nevertheless, in all the extremum seeking related work, only a single vehicle is used to collect measurements at different locations which is time-consuming. In addition, the trajectories generated by extremum seeking always demand costly maneuvers. Authors of \cite{Multi_deployment} propose strategies to deploy a group of vehicles around the source while applying extremum seeking. However, their focus is on achieving a formation distribution of the vehicles in accordance with signal strength instead of taking advantage of the vehicle swarm in collecting measurements. Therefore, this work is not a favorable example of multi-agent source seeking.

The advantages of robot swarm and sensor network attract many researchers to study them and apply them to the source seeking problem.  In \cite{CircularFormation, ConsensusCircularFormation} and \cite{MultiUAVMovingSource}, a team of agents implement a consensus algorithm to maintain a particular formation to track the gradient of the source. This method assumes the formation of the swarm being maintained perfectly which is too ideal to achieve. Cooperative source seeking algorithms are proposed in \cite{WZhang_Switching, Coop_control, veh_network} and \cite{Choi2012}. In \cite{WZhang_Switching}, the authors provide algorithms and experimental validation of a switching strategy for a team of agents trying to localize a source. Each agent switches from individual exploration to cooperative exploration only when individual gradient estimate is not available. A distributed coordination algorithm based on adaptive control is proposed in \cite{Choi2012}. However, it does not take noise into consideration which is crucial to our problem. In \cite{Coop_control} and \cite{veh_network}, cooperative control is applied to deploy agents in a way which is optimal for gradient estimation. But these two methods together with all other methods in this paragraph are all prone to be trapped in a local minimum. Since all of them deploy the swarm or network in a close neighborhood, they lack the global information of the model. Authors of \cite{pappas_stochastic} and \cite{Atanasov2014} apply stochastic approximation to the problem and enable the swarm to find the source in complex and noisy environment. But the computation complexity of the method hinders its implementation on some cheaper and less capable robots.
There are other source seeking methods that are non-gradient based. Some are developed by obtaining source functions. In \cite{Point_source, Source_obstacles, ElBadia2002, ElBadia2005} and \cite{Komornik2005}, researchers formulate the source seeking problem as an inverse problem. Depending on the source types, heat and wave equations are commonly used as candidate source functions. Parameters of source functions are found by optimizing the difference between collected data and simulated data from the candidate function. The source can be located after obtaining the source function. For instance, in \cite{Source_obstacles}, the incoming directions of waves are obtained after solving the inverse problem. By tracing back the wave direction rays, the source is located at the intersections of the rays. This method cannot be applied to our problem because it requires a priori knowledge about the candidate function that governs the decay profile of the source, and some information about the source like wave length and frequency. In our problem, only signal strength can be measured, and the signal decay profile is unknown. For the same reason, the statistical signal-processing technique `independent component analysis' \cite{Albini2003451} and statistical methods \cite{Chemical_plume, Optimotaxis}, \cite{Following_RF} used to construct a maximum likelihood map of the source location cannot be used for our problem. In \cite{IslerTrackingFish}, researchers use directional an tennas to obtain bearing measurements of tagged fish, and localize its location by triangulation. The method has been applied effectively in localizing the invaded carp in a lake. But it is not applicable to our problem since we can only access scalar measurements instead of bearing ones. In this paper, we address this problem using Particle Swarm Optimization (PSO) which does not require any a priori knowledge of the signal model emanating from the source.

PSO is a heuristic non-gradient based strategy, an evolutionary computation technique. It was first proposed by Kennedy and Eberhart \cite{PSO} who were inspired by the behavior of bird flock and fish school. Ever since then many variations of PSO have been proposed by researchers, like inertia weight PSO \cite{ModifiedPSO}, constriction PSO \cite{Clerc_Constriction}, neighborhood PSO \cite{NeighborPSO} in its early time, and Quantum behaved PSO \cite{QuantumPSO} and Digital Pheromone PSO \cite{Vijay_digital} developed recently.

As a swarm optimization technique, PSO has been applied to some source seeking tasks involving mobile robots. In \cite{PSOPugh2006} and \cite{PSOPugh2007}, PSO is modified to adapt to multi-robot search. The authors discuss the limitation posed by physical robots and conduct simulations with several communication models. However, these simulations are limited to some benchmark functions rather than real world signal sources. Since real sources such as electromagnetic signals, odor and heat sources are considerably different from the benchmark functions, the results in these papers are not remarkably useful to implement. In \cite{PSOJatmiko} and \cite{PSOMarques}, the authors incorporate potential field-based motion planning method with PSO, and propose strategies of localizing static and moving odor source in complex environment. While their strategies have been proven effective in simulations for localizing odor sources, we would like to focus our attention to electromagnetic sources. We will also explore different variations of PSO and compare the performance of various parameter configurations, topology models and obstacle avoidance strategies. In essence, our work focuses on finding the most effective PSO variation to solve the electromagnetic source seeking problem, and validate it in real experiments. Authors of \cite{PSODerr} also use RF signals as sources to be sought. But this work does not consider a more complex environment where obstacles exist, and is limited to simulations. A method for obtaining optimal PSO parameters is proposed in \cite{PSODoctor}. Essentially, it applies PSO at a lower level to seek for the source, and at a higher level to seek for the optimal parameter configuration. This iterative way of finding optimal parameters is reasonable computationally, but not applicable to real robots. It requires a significant number of trials at the lower level which leads to enormous amount of experimental data. The effort required to obtain the optimal parameter configuration in a specific scenario is extravagant for the simple goal of finding the source, especially when this optimal configuration can hardly be applied to a different scenario. \cite{HerefordPSO2007} provides some experimental results of implementing a modified PSO on real robots. They use a diffuse light source as the target and provide a simple strategy of dealing with obstacles. The experiments illustrate the efficacy of implementing PSO on real robots, but are constrained to a specific PSO configuration. In this work, we present extensive simulation and experimental results of various PSO variations and configurations, and provide suggestions on parameter selection.

This work is based on our previous exploration \cite{PSO_Rui_AIM}, \cite{PSO_Rui_MSC} on applying PSO to the source seeking problem. The main contributions of this paper are as follows. 1) We use a non-gradient based technique for the source seeking problem due to the inherent irregularity in the signal model. 2) We incorporate physical constraints posed by robots in the implementation and evaluation of PSO. 3) Guidelines are presented to choose proper parameters for several PSO variations. A strategy which enables PSO to be implemented experimentally in a complex environment is first presented in this paper.

This paper is outlined as follows. In Section \ref{sec:background}, we provide some background information regarding the problem description and the main concept of PSO. In Section \ref{sec:variations}, we evaluate and compare three PSO variations with different parameters. In Section \ref{sec:obstacles}, we propose collision avoidance strategies for implementing PSO in environments containing obstacles. In Section \ref{sec:experiments}, we present a description of the experimental setup, and discuss the implementation results. Finally, we conclude this work in Section \ref{sec:conclusion}.

\section{Background}
\label{sec:background}
\subsection{Problem Description}
Consider a point source located on a plane continuously transmits/emits a signal. Based on the assumption that a static source is present in the vicinity, a group of mobile agents, called {\it{seekers}}, explore the environment to locate the source. The scenario is similar to a colony of ant swarms trying to locate a food source. The seekers are assumed to be holonomic kinematic agents with maximum speed $v_{\max}$. The seekers have the capability to measure the strength of the signal emitted by the source at their current locations. However, the seekers have no information about the current location of the source, its signal strength and its decay profile. The objective of the seekers is to find the location of the source which is assumed at the location where the signal strength is maximum.

For most sources, signal intensity normally decays radially as the distance to the source increases. The decay profiles of some common sources are shown below:
\begin{itemize}
\item 
Let $P$ denote the power at which an electromagnetic source emits a signal. The decay profile of the signal intensity is given by the following equation \cite{rawnet11}
\begin{equation}
{P_A} = \frac{{cP}}{{(1 + d)}^\alpha },
\label{eqn:pow}
\end{equation}
where $P_{A}$ is the power of the signal measured at a point $A$ on the plane located at a distance $d$ from the source, $c$ and $\alpha$ are constants that depend on the physical parameters of the medium through which the signal is transmitted.

\item 
The concentration of a chemical $c(\vec{r}, t)$ at a point $\vec{r}$ emitted from a point source located at $\vec{\rho}$ emitting vapors at a constant rate of $\mu$ Kg/s is given by \cite{VaporSource}
\begin{equation}
\label{eqn:vapor}
c\left( {\vec r,t} \right) = \frac{\mu }{{4\pi \kappa \left| {\vec r - \vec \rho  } \right|}}{\text{erfc}}\left( {\frac{{\left| {\vec r - \vec \rho  } \right|}}{{2\sqrt {\kappa (t - {t_0})} }}} \right),
\end{equation}
where $\kappa$ is constant diffusivity in $m^2/s$. If we ignore the complementary error function ${\text{erfc}}(x) = (2/\sqrt \pi  )\int_x^\infty  {{e^{ - {y^2}}}dy}$ and only consider the dominant part on its left in equation (\ref{eqn:vapor}), the substance concentration is inversely proportional to the distance between $\vec{r}$ and the vapor source.

\item 
For a spherical sound source, the acoustic intensity $I_r$ at a point with a distance of $r$ from the center in the radial direction is given by \cite{PASPWEB2010}
\begin{equation*}
I_r = \frac{P}{4\pi r^2}
\end{equation*}

\end{itemize}

However, in reality the measured signal intensity is too noisy to be accurately described by these decay profiles. For instance, reflection, refraction, multi-path fading, etc. can influence the decay profile dramatically making the actual one highly different from the theoretical one.

\begin{figure}[htb]
	\centering
	\begin{subfigure}[b]{0.9\linewidth}
		\includegraphics[width=\textwidth]{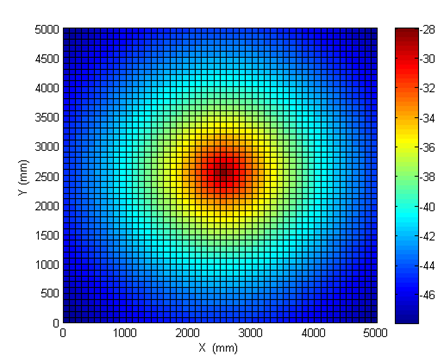}
	\end{subfigure}\\%
	\begin{subfigure}[b]{0.9\linewidth}
		\includegraphics[width=\textwidth]{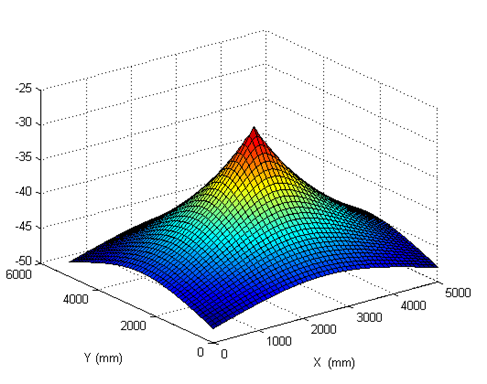}
	\end{subfigure}
	\caption{Theoretical decay profile of an electromagnetic source}
	\label{fig:Map_theory}
\end{figure}

Figure \ref{fig:Map} illustrates a real RSSI (Received Signal Strength Indication) map of an RF source provided by an XBee\textregistered  ZB RF module on a 5 m$\times$5 m plane. The XBee module was located at the center of this area, and the measurements were taken by another XBee module. The figure clearly illustrates the fact that the real RSSI profile has many local extrema and is non-differentiable almost everywhere contrary to the theoretical decay profile described by (\ref{eqn:pow}) shown in Figure \ref{fig:Map_theory}.

Therefore, an optimization method not limited by differentiability requirement, yet able to search highly multi-modal design spaces is desired for direct RSSI measurements, as portrayed in Figure 1. Also since mobile robots are used within the environment to seek maximum signal source, a population based method where each population member has a  one to one correspondence with a mobile robot is favored. Particle Swarm Optimization (PSO), as described in Section II B, has the ingredients required to address the above challenges.

\subsection{Original PSO}
\label{sec:Original}
In this subsection, we provide a brief description of the concept of PSO. PSO is a population based search algorithm first proposed in \cite{PSO} by Kennedy and Eberhart through simulation of a simplified social model. Although PSO was originally designed to solve minimization problems, it can be used to find the maximum of a function, with a simple change. It is initialized with a number of random solutions, called \textit{particles}. Each particle is also randomly initialized with a velocity within some user designed range. Each particle evolves iteratively in the search-space trying to improve the solution in the following manner:
\begin{eqnarray}
v^k_{i+1} &=& v^k_i +  U(0,c_1) (Pbest^k - x^k_i) \nonumber \\
          & &+ U(0,c_2) (Gbest - x^k_i) \label{eqn:PSO_basic}\\
x^k_{i+1} &=& x^k_i + v^k_{i+1} \nonumber \\
\end{eqnarray}
where $x^k_{i+1}$ and $v^k_{i+1}$ represent the position and velocity of the $k$th particle in the $i+1$th iteration, $U(0,c_1)$ and $U(0,c_2)$ are uniformly distributed random numbers within $[0, c_1]$ and $[0, c_2]$, and $Pbest$ and $Gbest$ are best previous position of a particle and best previous position in the swarm. A best previous position is where a particle obtains the minimum cost in its search history. In our case, the above equations can be interpreted in the following way: Assuming $n$ seekers as $n$ particles moving in the search-space $X$, the position of the $k$th seeker in the $i$th iteration is denoted as $x^k_i \in X \subset \mathbb{R}^2$. The cost function $f: \mathbb{R}^2 \rightarrow \mathbb{R}$ incurred by each seeker is the negative of the signal strength received at its current location. The objective of the seekers is to communicate, and move in a manner so as to reach the global minimum of the cost function.

We initialize the position and velocity of each seeker with a uniformly distributed vector $x^k_1$ in the search-space and $v^k_1$ within given bound. Each seeker is assumed to have the knowledge of its own best previous position and the global best previous position based on the assumption that each seeker has the memory to store its own previous experience and can benefit from the previous experience of all other members. Therefore, in (\ref{eqn:PSO_basic}), velocity $v^k_{i+1}$ consists of three terms: the effect of seeker's previous velocity, its best known position and global best known position. 

\section{PSO Modifications and Variations}
\label{sec:variations}
In this section, we will introduce some physical constraints into PSO and compare the performance of three different PSO variations and provide guidelines on parameter selection.

First, we shall consider the bounds to the search space. As an optimization technique, boundary conditions exist in PSO. However, different actions are taken when particles violate boundary conditions: some discard these particles, some bounce them back, and some confine them to the boundary, etc. Discarding seekers may impair the performance especially when the total number of seekers is limited and every seeker is of significant value to performance. Therefore, we choose to confine seekers to boundaries, namely,

Constraint 1: If $x^k_{i+1} \notin X$, then  $x^k_{i+1}$ is set to the boundary point on $X$ in the direction of $v^k_{i+1}$.

In PSO, there's no constraint on the velocity of a particle. It is possible for a particle to fly across the entire search space in a single iteration. However, this does not apply to seekers in our case which are actually ground robots. It is more appropriate to treat the velocity of a particle as a step of a robot in our implementation which decomposes the step length into speed and duration. Given a sufficiently long duration, a robot also can move a large step which crosses the entire search space. However, it is at the expense of a longer searching time, and a greater energy consumption, which is crucial to a robot with limited battery capacity. On the contrary, the step length should not be too small to affect the performance of PSO. To balance performance and efficiency, simulations to find a proper step length are conducted below. We denote the step length by $v_{\max}$, and check the following constraint in every iteration.

Constraint 2: If $|v^k_{i+1}| > v_{\max}$, then $|v^k_{i+1}| = v_{\max}$ with the direction of $v^k_{i+1}$ unchanged.

The aforementioned constraints apply to all simulations and experiments in this paper.

\subsection{PSO with Inertia Weight}
\label{sec:InertiaWeight}
One variation of the original PSO is to introduce an inertia weight $\omega$ to the previous velocity in (\ref{eqn:PSO_basic}), which leads to the following equation to update velocity \cite{ModifiedPSO},
\begin{equation}
v^k_{i+1} = \omega_i v^k_i + U(0,c_1) (Pbest^k - x^k_i)+ U(0,c_2) (Gbest - x^k_i)
\label{eqn:PSO_v}
\end{equation}
          
According to Shi and Eberhart's analysis in \cite{ModifiedPSO}, the inertia weight is critical in  balancing global and local search. If $\omega$ is set to zero, the seekers become ``memoryless" about its past velocity. With seekers' velocity only determined by individual and global best previous positions, all seekers would converge to the global best position directly making the search process resemble a local search. On the contrary, if $\omega$ is set to a larger number, the seekers are more stubborn in their previous velocity, which leads them to larger area of exploration. In other words, a larger inertia weight facilitates global exploration while a smaller one facilitates local exploitation to fine-tune the current search area \cite{PSOPara}. Therefore, implementing a damping mechanism to $\omega$ contributes to better global exploration in the beginning stage and better local exploitation when the swarm is closer to the source. 

The study in \cite{kennedy2001swarm} shows that $c_1$ and $c_2$ together contribute to the oscillation behavior for the seekers. As the values of $c_1$ and $c_2$ are increased, the frequency of oscillation of the seekers' trajectories also increases. Hereafter, we set $c_1 = c_2 = 2$ as suggested in \cite{kennedy2001swarm}. We start with multiplying $\omega$ with a damping coefficient $\lambda_{\omega}$ as the damping mechanism, and set $\lambda_{\omega} = 0.95$ as suggested in \cite{Vijay_digital}. (\ref{eqn:omega}) is implemented in every iteration after velocity is updates. Therefore,
\begin{equation}
\label{eqn:omega}
\omega_{i+1}  = \lambda_{\omega}\omega_i, \quad \text{with} \quad \lambda_{\omega} = 0.95
\end{equation}
%

\begin{figure}[thbp]
	\centering
	\begin{subfigure}{0.9\linewidth}
		\includegraphics[width=\textwidth]{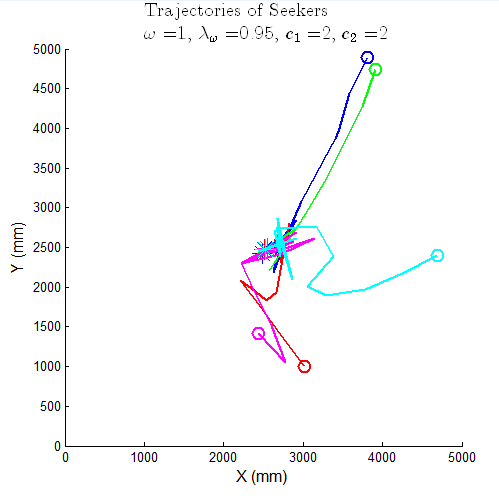}
		\caption{Trajectories of seekers with $\omega = 1$}
	\end{subfigure}\\%
	\begin{subfigure}{0.9\linewidth}
		\includegraphics[width=\textwidth]{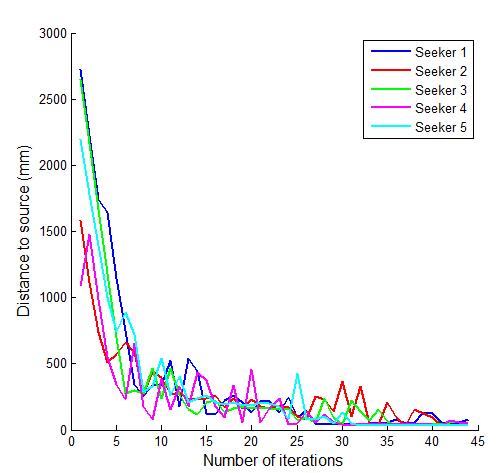}
		\caption{Statistics with $\omega = 1$}
	\end{subfigure}
	\caption{Trajectories of seekers with $\omega = 1$}
	\label{fig:traj_w1}
\end{figure}

\begin{figure}[thbp]
	\centering
	\begin{subfigure}{0.9\linewidth}
		\includegraphics[width=\textwidth]{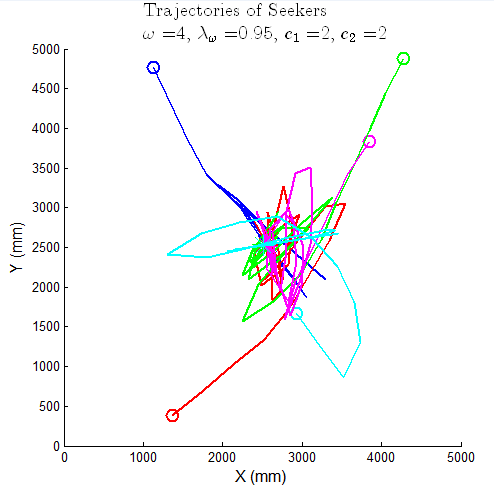}
		\caption{Trajectories of seekers with $\omega = 4$}
	\end{subfigure}
	\begin{subfigure}{0.9\linewidth}
		\includegraphics[width=\textwidth]{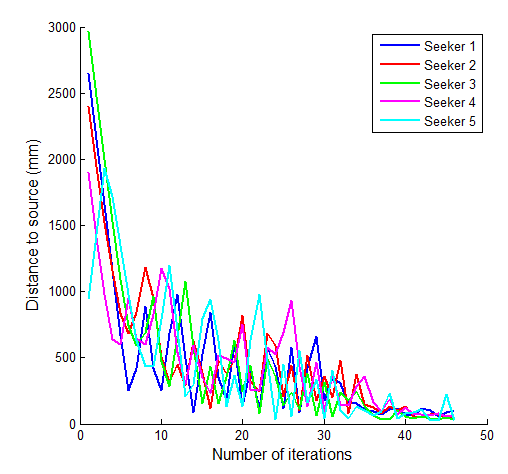}
		\caption{Statistics with $\omega = 1$}
	\end{subfigure}\\%
	\caption{Trajectories of seekers with $\omega = 4$}
	\label{fig:traj_w4}
\end{figure}

We choose the swarm size to be five in Section \ref{sec:InertiaWeight} and \ref{sec:Constriction}, and this will be explained in Section \ref{sec:SPSO}. Six sets of simulations with different initial $\omega$ and $v_{\max}$ were conducted. In each set, we ran 1000 simulations on the real RSSI design space described by Figure \ref{fig:Map}. The cost function is defined as the negative RSSI at each point which needs to be minimized. Each simulation was terminated when $Gbest$ remained unchanged for 20 iterations. Since the signal strength at the source is -28 dBm, we compared $Gbest$ with 28 after each simulation. In addition, we counted the number of iterations $I$ and the total distance traveled by all robots $TotalD$. Additionally, the following data is also collected:
\begin{itemize}
	\item $avgGbest$ denotes the mean of $Gbest$.
	\item $stdGbest$ denotes the standard deviation of $Gbest$.
	\item $avgI$ denotes the mean of $I$. 
	\item $avgTotalD$ denotes the mean of $TotalD$.
\end{itemize}
Simulation results are shown in Table \ref{table:lamda}, where the units of $v_{\max}$, and $avgTotalD$ are mm/iteration and mm, respectively.
\begin{table}[thbp]
	\centering
	\caption{Simulation results with different $\omega$ and $v_{\max}$, and a damping coefficient $\lambda_{\omega} = 0.95$ }
	\label{table:lamda}
\begin{tabular}{|c|c|c|c|c|c|c|}
\hline Set & $\omega_1$ & $v_{\max}$ & $avgGbest$ & $stdGbest$ & $avgI$ & $avgTotalD$ \\ 
\hline 1 & 2 & 500 & 28.1095 & 0.8095 & 45.639 & 52017 \\ 
\hline 2 & 3 & 500 & 28.0584 & 0.5715 & 46.613 & 66716 \\ 
\hline 3 & 4 & 500 & 28.0358 & 0.3631 & 47.187 & 77312 \\
\hline 4 & 5 & 500 & 28.0215 & 0.1187 & 47.855 & 84749 \\  
\hline 5 & 2 & 1000 & 28.0378 & 0.3553 & 48.313 & 101966 \\ 
\hline 6 & 3 & 1000 & 28.0272 & 0.1507 & 50.747 & 133442 \\ 
\hline 
\end{tabular} 
\end{table}
The first four sets illustrate the effect of increasing $\omega_1$ on the performance of the searching algorithm. As $\omega_1$ increases from 2 to 5, $avgGbest$ gets closer to 28, which means the seekers perform better in locating the source. Meanwhile, a decreasing $stdGbest$ represents growing reliability of the algorithm which is another indicator of improved performance. This improvement can be supported by the fact that larger $\omega$ facilitated global exploration. With a larger initial $\omega$, seekers tend to preserve their previous velocity, and explore a larger area in early iterations. Therefore, they are less likely to be trapped in a local minimum, and more likely to find the global minimum. However, the improved performance is at the expense of higher energy consumption. Though the average iterations $avgI$ is not clearly related to the change of $\omega_1$, $avgTotalD$ in set 4 is about 2.5 times that of set 1. Figures \ref{fig:traj_w1} and \ref{fig:traj_w4} are good demonstrations of the reason, in which seekers are represented by different colors and the small circles represent initial positions. It is clear from the figures that trajectories with $\omega_1=4$ cycle in a larger area and converge slower than that of $\omega_1 = 1$. So we can clearly see a trade-off between performance and energy consumption.

Sets 5 and 6 are used to reveal the influence of $v_{\max}$. Comparing sets 5 and 6 with 1 and 2, we find slight improvement in the $avgGbest$ and $stdGbest$ when $v_{\max}$ doubles. However, the average total distance traveled also doubles. Moreover, if we take set 4 into consideration, apparently, increasing $\omega$ is a better strategy than increasing $v_{\max}$ in terms of both performance and energy efficiency.

\subsection{PSO with Constriction Factor}
\label{sec:Constriction}

Another variation similar to PSO with inertia weight that is implemented to the source seeking problem in this paper is PSO with a constriction factor. Introduced by Clerc in \cite{Clerc_Constriction}, the constriction factor is used to prevent ``explosion" and ensure convergence of PSO. Equation (\ref{eqn:PSO_constriction}) and (\ref{eqn:K_constriction}) describe the basic concept of the constriction factor. 
\begin{equation}
\label{eqn:PSO_constriction}
v^k_{i+1} = K[v^k_i + U(0,c_1) (Pbest^k - x^k_i)+ U(0,c_2) (Gbest - x^k_i)]
\end{equation}
\begin{equation}
\label{eqn:K_constriction}
K = \frac{2}{\left|2-\phi - \sqrt{\phi^2 - 4\phi} \right| },\quad \text{where} \; \phi=c_1+c_2,\; \phi >4
\end{equation}
Compared to the original PSO, the entire RHS of (\ref{eqn:PSO_basic}) is multiplied by a coefficient $K$, called the constriction factor. $K$ is a function of $c_1$ and $c_2$ as shown in (\ref{eqn:K_constriction}). The main idea of constriction PSO is to take advantage of the mathematical nature of (\ref{eqn:K_constriction}) which guarantees the convergence of the algorithm. Detailed explanation of the mechanism of constriction PSO can be found in \cite{Clerc_Constriction} which is beyond the scope of this paper.

A closer look at (\ref{eqn:PSO_constriction}) reveals that it is a special case of (\ref{eqn:PSO_v}), whose inertia weight $\omega$ is set to $K$, and $c_1$ and $c_2$ are multiplied by $K$. It is the relation between $\phi$ and $K$ that prevents the swarm from ``explosion". Therefore, according to Clerc, $v_{\max}$ is not necessary when the constriction factor is applied. However,
for application reasons that were mentioned before, we would keep $v_{\max}$ to a smaller value to improve the energy efficiency of the robots. Since $K$ is a decreasing function of $\phi$ whose supremum is 1. If we think of it in terms of (\ref{eqn:PSO_v}), the supremum of $\omega$ is 1. This suggests the constriction PSO does not emphasize global exploration at the initial stage of the search. And it does not favor local exploitation, either, since $K$ does not vary through the search.

We also conducted 6 sets of simulations on the constriction PSO algorithm. In all sets, $c_1$ and $c_2$ are set to the same value of $\phi/2$ to balance the influence of individual and swarm experience. All configurations are identical to those in the previous section unless otherwise specified. Data is collected in Table \ref{table:constriction}.
\begin{table}[thbp]
	\centering
	\caption{Simulation results of constriction PSO}
	\label{table:constriction}
\begin{tabular}{|c|c|c|c|c|c|c|c|}
	\hline Set & $K$ & $\phi$ & $v_{\max}$ & $avgGbest$ & $stdGbest$ & $avgI$ & $avgTotalD$ \\ 
	\hline 7 & 0.5 & 4.5 & 500 & 29.2826 & 2.4494 & 51.443 & 15837 \\ 
	\hline 8 & 0.73 & 4.1 & 500 & 28.6432 & 1.9508 & 54.338 & 32282 \\ 
	\hline 9 & 0.8 & 4.05 & 500 & 28.5366 & 1.6601 & 49.583 & 44542 \\ 
	\hline 10 & 0.90 & 4.01 & 500 & 28.3358 & 1.4592 & 43.043 & 62037 \\ 
	\hline 11 & 0.73 & 4.1 & 1000 & 28.2765 & 1.1675 & 48.535 & 52972 \\ 
	\hline 12 & 0.8 & 4.05 & 1000 & 28.1730 & 0.9260 & 45.996 & 77154 \\ 
	\hline 
\end{tabular} 
\end{table}
Sets 7, 8, 9 and 10 show the impact of decreasing $\phi$, or increasing $K$. As $K$ increases, growing emphasis is put to the term of previous velocity in (\ref{eqn:PSO_constriction}). Therefore, the seekers tend to explore larger area and have a higher chance of finding the source. This improvement in performance is evident in these 4 sets, as both $avgGbest$ and $stdGbest$ decrease with $K$. And we can see $avgTotalD$ also grows with $K$ regardless of how $avgI$ varies. This is also the result of favoring global exploration, since the seekers ``fly" longer distance in each iteration when they emphasize exploration.

In set 11 and 12, we keep $K$ equal to that in set 12 and 13, and only double $v_{\max}$. We can see significant improvement in performance when $v_{\max}$ doubles. Because this allows the entire RHS of (\ref{eqn:PSO_constriction}) to be doubled, including the term for previous velocity which puts emphasis on global exploration in another way. However, this improvement is not seen in set 5 and 6 in PSO with inertia weight. A reasonable guess may be that when the coefficient $\omega$ is large, the performance is mainly influenced by $\omega$ rather than $v_{\max}$.

\begin{figure}[thbp]
	\centering
	\begin{subfigure}{0.9\linewidth}
		\includegraphics[width=\textwidth]{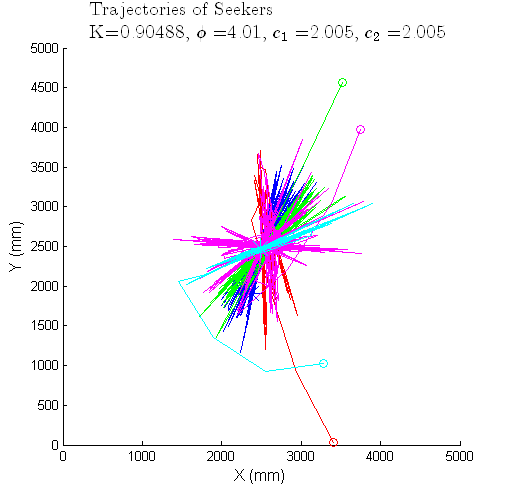}
		\caption{Trajectories}
	\end{subfigure}
	\begin{subfigure}{0.9\linewidth}
		\includegraphics[width=\textwidth]{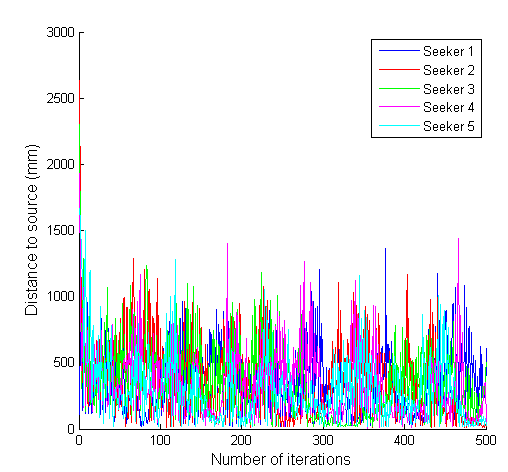}
		\caption{Statistics}
	\end{subfigure}%
	\caption{Trajectories of seekers with a constriction factor of $K = 0.9$}
	\label{fig:pso_oscillation}
\end{figure}

Observation on the trajectories of seekers reveals another feature of constriction PSO which can be seen in Figure \ref{fig:pso_oscillation}. This figure illustrates one typical simulation result where $K = 0.9$. In this simulation, the source is found after about 40 iterations, however, the swarm does not converge to the source after that as in Figures \ref{fig:traj_w1} and \ref{fig:traj_w4}. Instead, all seekers keep oscillating around the source showing no sign of convergence. Comparing this to those sets on PSO with inertia weight, we can find that the violent oscillation in Figure \ref{fig:pso_oscillation} actually roots in the lack of a damping mechanism in constriction PSO. With a constant coefficient for the previous velocity, the swarm is incapable of switching from favoring global exploration in initial stage of the search to favoring local exploitation in later stage.

Based on the collected the data, our preliminary judgment is PSO with inertia weight is better suited for our application. Moreover, 1/10th of the length of the search space is a reasonable value for $v_{\max}$. As for the inertia weight, any value between 2 and 4 should produce some good results.

\subsection{SPSO}
\label{sec:SPSO}
The last PSO variation studied in this paper is Standard Particle Swarm Optimization (SPSO). It is a substantial improvement to the original PSO published in 1995, and researchers that developed their own PSO implementations benchmark their method’s performance against SPSO. The implementation of SPSO 2006 can be found here \cite{SPSO}. In this section, we will first provide a brief description of SPSO 2006, then study three SPSO topology models.

The velocity update equation in SPSO is almost the same as Equation (\ref{eqn:PSO_v}), except that $Gbest$ is replaced with $Lbest$ --  best previous position in the neighborhood, as shown in the following equation.
\begin{equation}
v^k_{i+1} = \omega v^k_i + U(0,c) (Pbest^k - x^k_i) + U(0,c) (Lbest^k - x^k_i) \label{eqn:SPSO}
\end{equation}
As a benchmark variation, there are generally accepted values for all the parameters in SPSO. The swarm size is determined by $10+[2\sqrt{D}]$, where $D$ is the dimension of the search space. So we use 12 seekers in this subsection. Other parameter values are 
\begin{eqnarray}
\omega &=& \frac{1}{2\ln(2)}\approx 0.721 \nonumber\\
c &=& \frac{1}{2}+\ln(2) \approx 1.193 \nonumber
\end{eqnarray}
Please refer to \cite{SPSO2011} for detailed description on initialization and confinement of SPSO.

A noticeable distinction in SPSO is the introduction of neighborhood. Neighborhood defines the communication topology among seekers. In this subsection, we will study and compare the implementation of three commonly used models on the source seeking problem.

\begin{figure}[htbp]
	\centering
	\begin{subfigure}[b]{0.28\textwidth}
		\includegraphics[width=\textwidth]{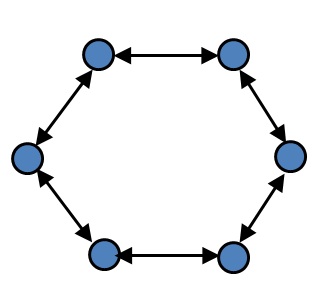}
		\caption{Ring topology}
		\label{fig:graph_ring}
	\end{subfigure}\\%
	~ 
	\begin{subfigure}[b]{0.3\textwidth}
		\includegraphics[width=\textwidth]{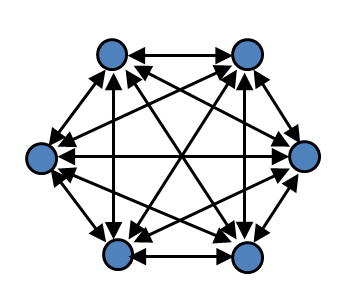}
		\caption{Fully connected topology}
		\label{fig:graph_full}
	\end{subfigure}\\
	~ 
	\begin{subfigure}[b]{0.32\textwidth}
		\includegraphics[width=\textwidth]{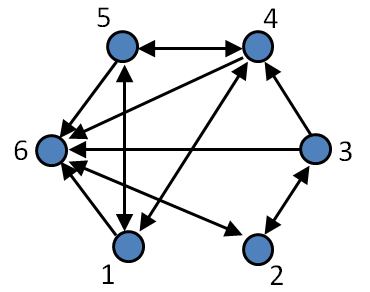}
		\caption{Adaptive random topology}
		\label{fig:graph_random}
	\end{subfigure}
	\caption{Graphs of different topologies}\label{fig:graph}
\end{figure}

Figure \ref{fig:graph} present the graphs of all three models. Figures \ref{fig:graph}(a) and \ref{fig:graph}(b) are self-explanatory. Figure \ref{fig:graph}(c) is the adaptive random topology model \cite{clerc2010particle} when $K = 3$. In this model, each particle informs $K$ random particles and itself of its $Pbest$, which means it informs at most $K+1$ different particles and at least one particle (itself). For instance, in Figure \ref{fig:graph_random}, particle 6 informs particle 2 and itself and has 5 informants $\{1, 2, 3, 4, 5\}$. $Lbest$ of a particle is defined as the best $Pbest$ among all its informants. This graph changes after every unsuccessful iteration (no improvement in $Gbest$).

To compare the aforementioned topology models, we conducted five sets of simulations with different models. Table \ref{table:topology} collect all simulation data.
 \begin{table}[htbp]
 	\centering
 	\caption{Simulation results with different topology models }
 	\label{table:topology}
 	\begin{tabular}{|c|c|c|c|c|c|}
 		\hline Set & Topology & $avgGbest$ & $stdGbest$ & $avgI$ & $avgTotalD$ \\ 
 		\hline 13 & ring & 28.000 & 1.82E-05 & 29.970 & 68380 \\ 
 		\hline 14 & fully connected & 28.000 & 6.74E-04 & 29.331 & 65475 \\ 
 		\hline 15 & $K = 3$ & 28.000 & 2.20E-04 & 29.259 & 68860 \\
 		\hline 16 & $K = 6$ & 28.001 & 3.16E-02 & 28.671 & 66913 \\  
 		\hline 17 & $K = 12$ & 28.002 & 4.47E-02 & 29.127 & 65212 \\  
 		\hline 
 	\end{tabular}
 \end{table}
Surprisingly, there's no distinguishable difference among these various models either in terms of $Gbest$ or $avgTotalD$. Consequently, we cannot draw any solid conclusion on the superiority of one model over the others. One plausible reason for this inconclusive result may lie in the number of seekers. 12 seekers seems to be excessive for our implementation making the influence of topology model and other parameters negligible. For the same reason, we only used five seekers in previous subsections to distinguish influence of those parameters of interest. In future implementations, we would prefer the fully connected model for simplicity reason. 

\section{PSO in Complex Environment}
\label{sec:obstacles}
In previous implementations, the source seeking task is carried out in an ideal obstacle-free environment. However, in real-world, we have to cope with obstacles as well as collisions among seekers which can be modeled as collision avoidance in the presence of dynamic obstacles. Therefore, we decompose the obstacle avoidance problem into two stages to deal with static and dynamic obstacles, respectively.

\subsection{Static Obstacles}

Static obstacles are common in a search environment. Constructions and uneven terrain are all potential static obstacles for seekers. We will give a short description of two static obstacle avoidance strategies proposed in our previous work \cite{PSO_Rui_AIM}, \cite{PSO_Rui_MSC}. Then we will integrate them into SPSO and compare their performance in simulations.

\begin{figure}[htbp]
	\centering
	\includegraphics[width=0.9\linewidth]{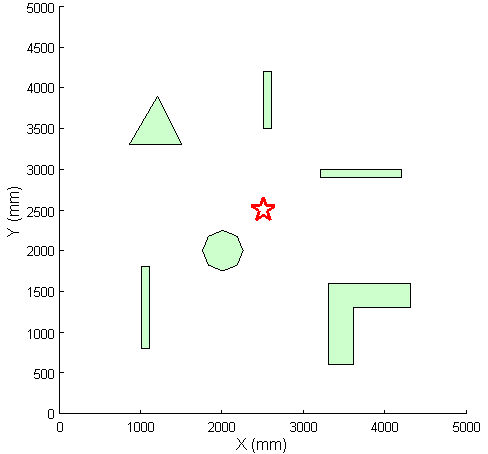}
	\caption{Map with obstacles}
	\label{fig:obstacle_map}
\end{figure}

\begin{figure}[htb]
	\centering
	\begin{subfigure}[b]{0.8\linewidth}
		\includegraphics[width=\textwidth]{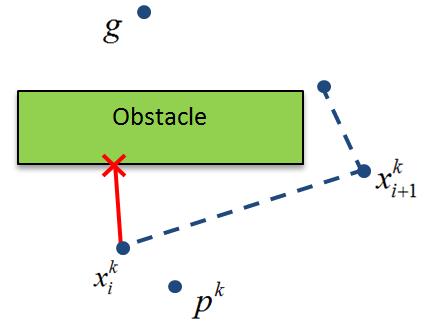}
		\caption{Strategy 1}
		\label{fig:random_step}
	\end{subfigure}\\%
	~ 
	\begin{subfigure}[b]{\linewidth}
		\includegraphics[width=\textwidth]{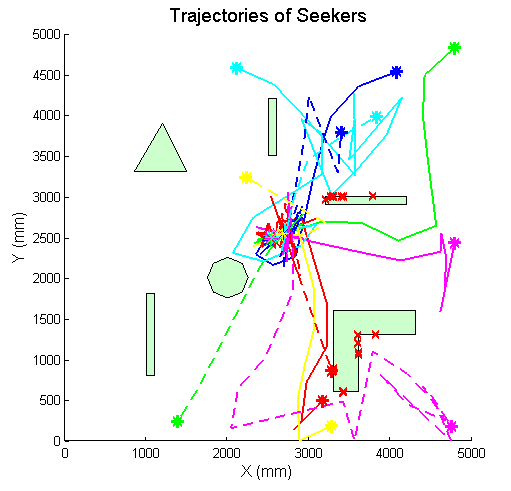}
		\caption{Trajectories of seekers}
		\label{fig:traj_rand}
	\end{subfigure}
	\caption{Static obstacle avoidance -- Strategy 1}
\end{figure}

Obstacles are described as simple convex or concave polygons in the search space as shown in Figure \ref{fig:obstacle_map}. The red star in the center represents the source. Seekers are provided with the information about each obstacle's position and size beforehand. The main idea of integrating obstacle avoidance into SPSO is to add a new operation mode to the seekers. They operate in the regular mode implementing SPSO when their trajectories do not collide with obstacles, and switch to obstacle avoiding mode when there's potential collision.
Strategy 1 inherits the heuristic nature of PSO. It introduces a step with a specific length and a random direction into PSO when an obstacle lies in the next step of a seeker. We set the length of this random step to be the ``diameter" of the obstacle so that the seeker has a good chance of circumvent the obstacle in one step as shown in Figure \ref{fig:random_step}. Here diameter refers the largest distance between any two points on the obstacle. Let $D_j$ denote the diameter of the $j$th obstacle. Algorithm \ref{algo:Strategy1} presents the procedure of this strategy. It is executed whenever a new step is generated for a seeker by PSO. In other words, collision with any obstacle is always checked for every step from $x^k_i$ to $x^k_{i+1}$ before it is executed. Figure \ref{fig:traj_rand} demonstrates the trajectories of 12 seekers implementing Strategy 1 in SPSO. Different seekers' trajectories are represented by different line styles. ``*" denotes the initial position of each seeker, and red ``x" represents a potential collision with an obstacle.
\begin{algorithm}[htbp]
	\caption{Static Obstacle Avoidance Strategy 1}
	\label{algo:Strategy1}
	\begin{algorithmic} [1]
		\IF {$x^k_{i+1}$ is in the $j$ the obstacle}
		\REPEAT
		\STATE set $v^k_{i+1}$ to a random direction and let $|v^k_{i+1}| = D_j$
		\STATE $x^k_{i+1} = x^k_i + v^k_{i+1}$
		\UNTIL {$x^k_{i+1}$ is not in any obstacle}
		\ENDIF
		\STATE Proceed with the normal PSO
	\end{algorithmic}
\end{algorithm}

\begin{figure}[htb]
	\centering
	\begin{subfigure}[b]{0.8\linewidth}
		\includegraphics[width=\textwidth]{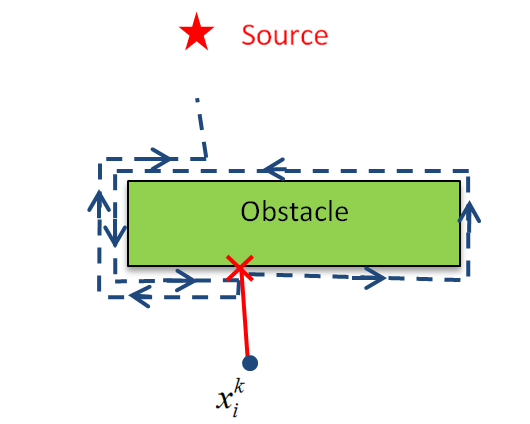}
		\caption{Strategy 2}
		\label{fig:bug1}
	\end{subfigure}\\%
	~ 
	\begin{subfigure}[b]{\linewidth}
		\includegraphics[width=\textwidth]{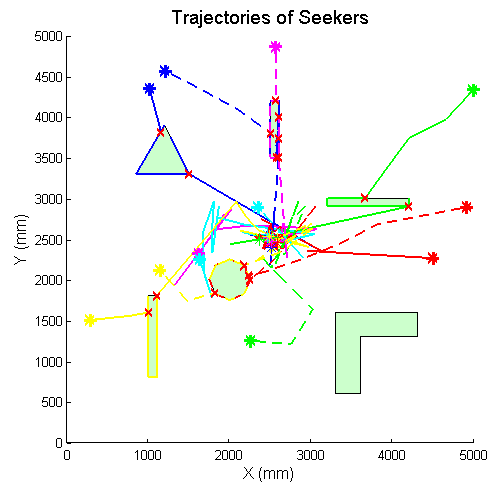}
		\caption{Trajectories of seekers}
		\label{fig:traj_bug}
	\end{subfigure}
	\caption{Static obstacle avoidance -- Strategy 1 (Bug 1 Algorithm)}
\end{figure}

Strategy 2 is a variation of the \textit{Bug 1} algorithm \cite{bug_algorithm}. Instead of knowing the position of the goal, only the signal strength at its current position is known to a seeker in our case. Once a seeker switches to obstacle avoidance mode, it starts to circumnavigate the encountered obstacle. As it circumnavigates, it measures the signal strength along its path. After circumnavigating the entire obstacle, the seeker follows the shortest path on the boundary to point at which it measures the largest signal strength and implements regular SPSO. Although in our case, it is not guaranteed that the seeker would end at the closest point to the source on the obstacle's boundary as in the \textit{Bug Problem}, it is highly likely to be on the side of the obstacle which is closer to the source. Because the source signal strength generally decays with distance, though it is quite noisy and does not strictly follow a decay profile. This provides the 
basis of implementing the \textit{Bug 1} algorithm and prevents the seeker from going back to the same obstacle. Figure \ref{fig:bug1} illustrates the trajectory of a seeker implementing the ``Bug 1" algorithm to avoid obstacle. And Figure \ref{fig:traj_bug} presents the trajectories of 12 seekers implementing Strategy 2 in SPSO.

Now we provide more simulation results to compare these two obstacle avoidance strategies. We conducted 4 sets of simulations. Set 18 and 19 used the parameters in set 2, and set 20 and 21 used the fully connected topology. Simulation results are collected in Table \ref{table:obstacle}.

\begin{table}[thbp]
	\centering
	\caption{Simulation results for two obstacle avoidance strategies}
	\label{table:obstacle}
	\begin{tabular}{|c|c|c|c|c|c|c|}
		\hline Set & Variation & Strategy & $avgGbest$ & $stdGbest$ & $avgI$ & $avgTotalD$ \\
		\hline 18 & Inertia & 1 & 28.0229 & 0.3779 & 45.3947 & 68475 \\ 
		\hline 19 & Inertia & 2 & 28.2634 & 1.2469 & 66.3155 & 66400 \\  
		\hline 20 & SPSO & 1 & 28 & 2.65E-04 & 33.776 & 74017 \\ 
		\hline 21 & SPSO & 2 & 28.0826 & 0.4034 & 35.867 & 50034 \\ 
		\hline 
	\end{tabular} 
\end{table}

Strategy 1 outperforms Strategy 2 in both $avgGbest$ and $stdGbest$ for both PSO variations. Very small standard deviation suggests the high reliability of Strategy 1. The reason Strategy 1 ends with longer distance is that its random step is usually larger than $v_{\max}$ because of the size of obstacles. While in Strategy 2, seekers usually take steps shorter than $v_{\max}$ when circumnavigating obstacles. Overall, Strategy 1 is better than Strategy 2 in simulations. 
Moreover, performance distinction originated from different variations is more significant than from different strategies. This primarily result from the size of the swarm.

\subsection{Dynamic Obstacles}

In all previous simulations, seekers are assumed to be points on a plane. However, in practice they have a finite area. This makes dynamic obstacles avoidance an inevitable issue in the application of swarm robots since every robot acts as a dynamic obstacle to others. To deal with this problem, we add two steps to the obstacle avoidance mode.

During each iteration, after $x^k_{i+1}$ are generated by PSO and checked or modified using the static obstacle avoidance strategy, potential collisions among seekers need to be checked. In this stage, there are two possible kinds of collisions: (1) Collisions at seekers end points; (2) Collisions in seekers trajectories. Since the seekers are assumed to dimensionless point particles in PSO, the algorithm needs to be modified to take into account possible collision between the robots at the end of their paths in a real scenario. Some seekers maybe too close to fit in the real robots causing collisions at these end points. In order to circumvent this problem, we incorporate a model that forces the seekers to repel each other to rearrange their end points to avoid collision. This is described in Algorithm \ref{algo:RepulsiveForce}.
\begin{algorithm}[htbp]
	\caption{End point arrangement using repulsive force}
	\label{algo:RepulsiveForce}
	\begin{algorithmic} [1]
		\STATE $S$ is the set of seekers
		\STATE $R$ is the radius of a seeker
		\STATE $t$ is a scaling factor
		\WHILE {$\exists \; |x^p_{i+1} - x^q_{i+1}| < 2R, \; p,q \in S, p\neq q$}
			\FOR{each $k \in S$}
					\FOR{each $j \in S, j \neq k$}
						\STATE $d = x^k_{i+1}-x^j_{i+1}$
						\IF{$|d| >= 2R$}
							\STATE $Force(k,j) = 0$
						\ELSE
							\STATE $Force(k,j)= d(2R-|d|)/|d|$
						\ENDIF
					\ENDFOR
				\STATE $Force(k) = \sum_{j \in S, j \neq k} Force(k,j)$
				\STATE $x^k_{i+1} = x^k_{i+1} + t Force(k)$
			\ENDFOR 
		\ENDWHILE
	\end{algorithmic}
\end{algorithm}

Algorithm \ref{algo:RepulsiveForce} ensures safe distance between any two seekers, and avoids end point collision. After this, if any seeker happens to lie in the path of others, the second step is activated. In this mode, seekers move sequentially. Only one seeker moves at a time while others stay still. We treat all other seeker as rectangular obstacles. We construct a reduced visibility graph \cite{choset2005principles} from the current position $x^k_i$ of the activated seeker to its next position $x^k_{i+1}$. Finally, by applying the Dijkstra's algorithm \cite{Dijkstra}, we generate the shortest path from $x^k_i$ to $x^k_{i+1}$.

\begin{figure}[htbp]
	\centering
	\includegraphics[width=\linewidth]{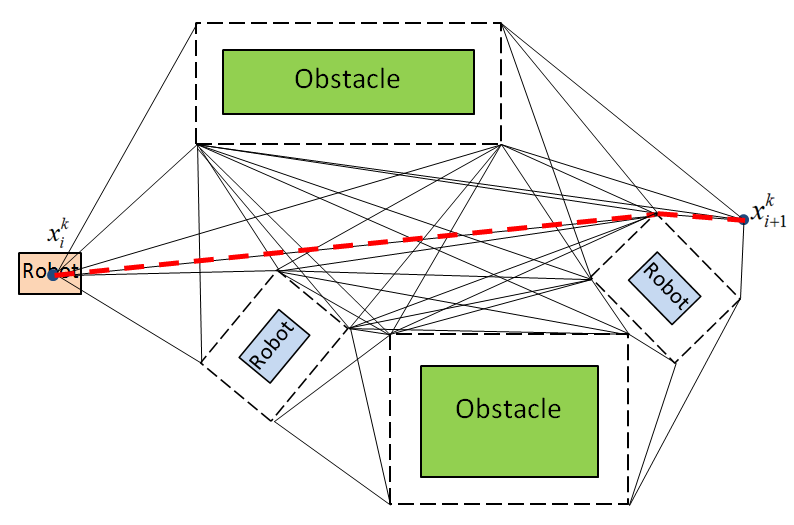}
	\caption{Visibility graph and shortest path}
	\label{fig:visibility}
\end{figure}

Figure \ref{fig:visibility} presents an example of the visibility graph and the shortest path. Due to the finite non-zero area of a seeker, the boundaries of obstacles and stationary seekers are expanded to the black dashed line to ensure a safety zone for the activated seeker (Minkowski sum of the obstacles with the seekers). The solid black lined delineate the visibility graph. The red dashed line represent the shortest path between $x^k_i$ and $x^k_{i+1}$.

So far, we have proposed a complete solution to implementing PSO on real robots in a complex environment where there exist potential collisions. In the next section, we will describe the experimental setup for implementation.

\section{Experiments}
\label{sec:experiments}

\begin{figure}[thbp]
 	\centering
 	\includegraphics[width=0.9\linewidth]{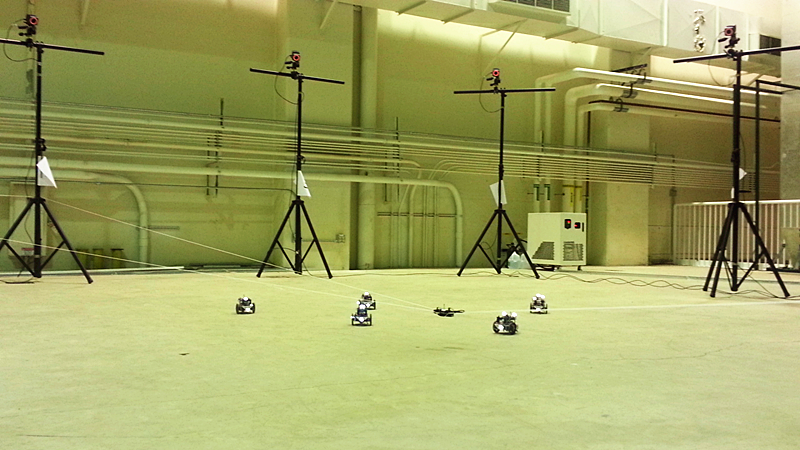}
 	\caption{Testbed}
 	\label{fig:setup}
 \end{figure}

\begin{figure}[thbp]
 	\centering
 	\includegraphics[width=0.9\linewidth]{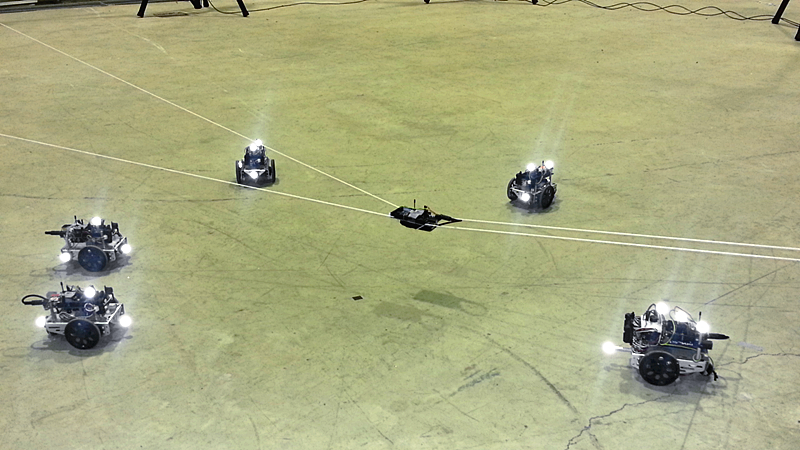}
 	\caption{Experiment environment with obstacles}
 	\label{fig:setup_close}
 \end{figure}

\begin{figure}[thbp]
  	\centering
  	\includegraphics[width=0.9\linewidth]{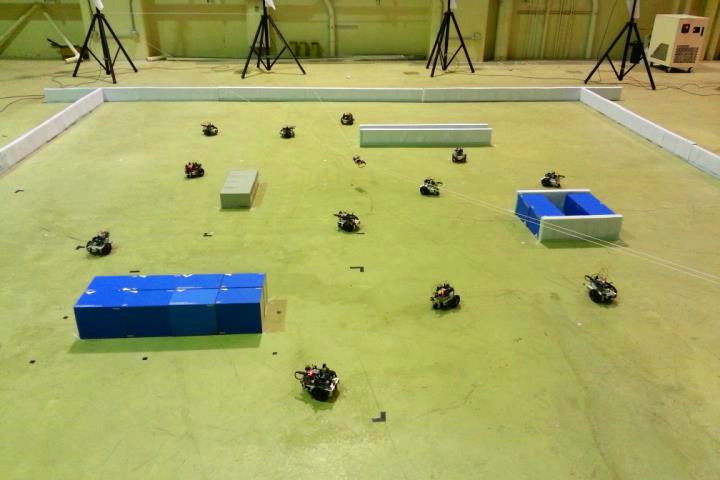}
  	\caption{Close look of robots and source}
  	\label{fig:complex_environment}
  \end{figure}

Our testbed is built on a 5 m$\times$5 m area covered by the Vicon tracking system. This system provides accurate position information of robots by recognizing markers on the robots serving as an indoor GPS system. The
source is an XBee module hanging in the middle at a height of 20 cm above floor. We do not place it on the floor in order to avoid potential collision with the robot. Robots used in experiments are small differential-drive robots modified from the Parallax Shield-Bot controlled by Arduino. Each robot is equipped with an XBee module to measure RSSI. Figures \ref{fig:setup} and \ref{fig:setup_close} are pictures of the testbed and robots. Figure \ref{fig:complex_environment} illustrates the complex environment with obstacles in which experiments were conducted.

 
 In the experiments, we built a centralized system with a computer being the center collecting and distributing information from and to all robots. This is not necessary since the strategies proposed in this paper are not computationally expensive and can be implemented on these robots without a strong computation ability. Also, since each robot can also communicate with each other, this system can work effectively without a central unit if robots have access to their own positions.
 
 Two successful experiments were recorded in the video. In these experiments, five robots were deployed to seek the source implementing the proposed strategies in an environment with obstacles. The parameters were chosen to be the same as set 2 in Table \ref{table:lamda}. 

\section{Conclusion}
\label{sec:conclusion}
In this paper, we explored the implementation of PSO to the electromagnetic source seeking problem. We modified PSO in accordance with the physical constraints posed by robots and the environment. Three PSO variations were evaluated through simulations. We found that the inertia weight PSO is best suited to our implementation and provided guidelines on parameter selection in PSO. We extended PSO from a pure computation technique to a complete solution to the source seeking problem in complex environment. Collision avoidance techniques were discussed extensively in this paper, and a complete obstacle avoidance strategy was incorporated in PSO. Our work was validated eventually in experiments using real robots.

In the future, we plan to explore and develop more advanced PSO variations that are specific for robotics applications. We would like to extend our work to more general source seeking scenarios, where sources may have different features and the obstacles in the environment cannot be simplified as polygons. Though it is unlikely that any variation can perform effectively in all kinds of scenarios, it is possible to explore the preferences of various scenarios and provide guidance in the selection of variations and parameter configurations.

%

%
%
%

%

\ifCLASSOPTIONcaptionsoff
  \newpage
\fi



\bibliographystyle{IEEEtran}
\bibliography{mybib}

\begin{thebibliography}{10}
\providecommand{\url}[1]{#1}
\csname url@samestyle\endcsname
\providecommand{\newblock}{\relax}
\providecommand{\bibinfo}[2]{#2}
\providecommand{\BIBentrySTDinterwordspacing}{\spaceskip=0pt\relax}
\providecommand{\BIBentryALTinterwordstretchfactor}{4}
\providecommand{\BIBentryALTinterwordspacing}{\spaceskip=\fontdimen2\font plus
\BIBentryALTinterwordstretchfactor\fontdimen3\font minus
  \fontdimen4\font\relax}
\providecommand{\BIBforeignlanguage}[2]{{%
\expandafter\ifx\csname l@#1\endcsname\relax
\typeout{** WARNING: IEEEtran.bst: No hyphenation pattern has been}%
\typeout{** loaded for the language `#1'. Using the pattern for}%
\typeout{** the default language instead.}%
\else
\language=\csname l@#1\endcsname
\fi
#2}}
\providecommand{\BIBdecl}{\relax}
\BIBdecl

\bibitem{Optimotaxis}
A.~R. Mesquita, J.~P. Hespanha, and K.~Astrom, ``Optimotaxis: A stochastic
  multi-agent optimization procedure with point measurements,'' in \emph{Hybrid
  Systems: Computation and Control, 11th International Workshop}, vol. 4981,
  2008, pp. 358--371.

\bibitem{AnimalNavigation}
N.~J. Vickers, ``Mechanisms of animal navigation in odor plumes,'' in
  \emph{Biol. Bull.}, vol. 198, no.~2, 2000.

\bibitem{TDoA}
B.~Charrow, N.~Michael, and V.~Kumar, ``Cooperative multi-robot estimation and
  control for radio source localization,'' in \emph{Experimental Robotics},
  ser. Springer Tracts in Advanced Robotics.\hskip 1em plus 0.5em minus
  0.4em\relax Springer International Publishing, 2013, vol.~88, pp. 337--351.

\bibitem{Coop_control}
P.~Ogren, E.~Fiorelli, and N.~Leonard, ``Cooperative control of mobile sensor
  networks:adaptive gradient climbing in a distributed environment,''
  \emph{IEEE Transactions on Automatic Control}, vol.~49, no.~8, pp.
  1292--1302, 2013.

\bibitem{Nehorai}
A.~Nehorai, B.~Porat, and E.~Paldi, ``Detection and localizing of
  vapor-emitting sources,'' \emph{IEEE Transactions on Signal Processing},
  vol.~43, no.~1, pp. 243--253, 1995.

\bibitem{VaporSource}
B.~Porat and A.~Nehorai, ``Localizing vapor-emitting sources by moving
  sensors,'' \emph{IEEE Transactions on Signal Processing}, vol.~44, no.~2, pp.
  1018--1021, 1996.

\bibitem{AUVCircle}
E.~Burian, D.~Yoerger, A.~Bradley, and H.~Singh, ``Gradient search with
  autonomous underwater vehicles using scalar measurements,'' in
  \emph{Proceedings of the 1996 Symposium on Autonomous Underwater Vehicle
  Technology}, 1996, pp. 86--98.

\bibitem{RotationBasedAngle}
S.~Venkateswaran, J.~Isaacs, K.~Fregene, R.~Ratmansky, B.~Sadler, J.~Hespanha,
  and U.~Madhow, ``Rf source-seeking by a micro aerial vehicle using
  rotation-based angle of arrival estimates,'' in \emph{American Control
  Conference (ACC)}, 2013, pp. 2581--2587.

\bibitem{SingleWingMAV}
K.~Fregene and C.~Bolden, ``Dynamics and control of a biomimetic single-wing
  nano air vehicle,'' in \emph{American Control Conference (ACC)}, 2010, pp.
  51--56.

\bibitem{Extremum_seeking}
K.~B. Ariyur and M.~Krstic, \emph{Real-Time Optimization by Extremum-Seeking
  Control}.\hskip 1em plus 0.5em minus 0.4em\relax Wiley-Interscience, 2003.

\bibitem{Non-h_forward_v}
C.~Zhang, D.~Arnold, N.~Ghods, A.~Siranosian, and M.~Krstic, ``Source seeking
  with nonholonomic unicycle without position measurement---part i: Tuning of
  forward velocity,'' in \emph{45th IEEE Conference on Decision and Control},
  Dec 2006, pp. 3040--3045.

\bibitem{Cochran2009}
J.~Cochran and M.~Krstic, ``Nonholonomic source seeking with tuning of angular
  velocity,'' \emph{Automatic Control, IEEE Transactions on}, vol.~54, no.~4,
  pp. 717--731, April 2009.

\bibitem{3-D_seeking}
J.~Cochran, A.~Siranosian, N.~Ghods, and M.~Krstic, ``3-d source seeking for
  underactuated vehicles without position measurement,'' \emph{IEEE
  Transactions on Robotics}, vol.~25, no.~1, pp. 117--129, 2009.

\bibitem{Liu20101443}
S.-J. Liu and M.~Krstic, ``Stochastic source seeking for nonholonomic
  unicycle,'' \emph{Automatica}, vol.~46, no.~9, pp. 1443 -- 1453, 2010.

\bibitem{Stanković20101243}
M.~S. Stanković and D.~M. Stipanović, ``Extremum seeking under stochastic
  noise and applications to mobile sensors,'' \emph{Automatica}, vol.~46,
  no.~8, pp. 1243 -- 1251, 2010.

\bibitem{Multi_deployment}
N.~Ghods and M.~Krstic, ``Multiagent deployment over a source,'' \emph{IEEE
  Transactions on Control Systems Technology}, vol.~20, no.~1, pp. 277--285,
  2012.

\bibitem{CircularFormation}
B.~Moore and C.~Canudas-de\mbox{-}Wit, ``Source seeking via collaborative
  measurements by a circular formation of agents,'' in \emph{American Control
  Conference (ACC)}, 2010, pp. 6417--6422.

\bibitem{ConsensusCircularFormation}
L.~Brinon-Arranz and L.~Schenato, ``Consensus-based source-seeking with a
  circular formation of agents,'' in \emph{European Control Conference (ECC)},
  2013, pp. 2831--2836.

\bibitem{MultiUAVMovingSource}
S.~Zhu, D.~Wang, and C.~Low, \emph{\BIBforeignlanguage{English}{Journal of
  Intelligent \& Robotic Systems}}, vol.~74, no. 1-2, pp. 333--346, 2014.

\bibitem{WZhang_Switching}
W.~Wu and F.~Zhang, ``Robust cooperative exploration with a switching
  strategy,'' \emph{IEEE Transactions on Robotics}, vol.~28, no.~4, pp.
  828--839, 2012.

\bibitem{veh_network}
R.~Bachmayer and N.~Leonard, ``Vehicle networks for gradient descent in a
  sampled environment,'' in \emph{Proceedings of the 41st IEEE Conference on
  Decision and Control}, vol.~1, 2002, pp. 112--117 vol.1.

\bibitem{Choi2012}
M.~Jadaliha, J.~Lee, and J.~Choi, ``Adaptive control of multiagent systems for
  finding peaks of uncertain static fields,'' \emph{Journal of Dynamic Systems,
  Measurement, and Control}, vol. 134, no.~5, 2012.

\bibitem{pappas_stochastic}
S.~Azuma, M.~S. Sakar, and G.~Pappas, ``Stochastic source seeking by mobile
  robots,'' \emph{IEEE Transactions on Automatic Control}, vol.~57, no.~9, pp.
  2308--2321, 2012.

\bibitem{Atanasov2014}
N.~A. Atanasov, J.~Le~Ny, and G.~J. Pappas, ``Distributed algorithms for
  stochastic source seeking with mobile robot,'' \emph{Journal of Dynamic
  Systems, Measurement, and Control}, 2014.

\bibitem{Point_source}
M.~Burger, Y.~Landa, N.~Tanushev, and R.~Tsai, ``Discovering a point source in
  unknown environments,'' in \emph{Algorithmic Foundation of Robotics VIII},
  ser. Springer Tracts in Advanced Robotics.\hskip 1em plus 0.5em minus
  0.4em\relax Springer Berlin Heidelberg, 2009, vol.~57, pp. 663--678.

\bibitem{Source_obstacles}
Y.~Landa, N.~Tanushev, and R.~Tsai, ``Discovery of point sources in the
  helmholtz equation posed in unknown domains with obstacles,'' \emph{Comm. in
  Math. Sci.}, vol.~9, no.~3, pp. 903--928, 2011.

\bibitem{ElBadia2002}
A.~E. Badia and T.~Ha-Duong, ``On an inverse source problem for the heat
  equation. application to a pollution detection problem,'' \emph{Journal of
  Inverse and Ill-Posed Problems}, vol.~10, no.~6, pp. 585--599, 2002.

\bibitem{ElBadia2005}
A.~E. Badia, T.~Ha-Duong, and A.~Hamdi, ``Identification of a point source in a
  linear advection–dispersion–reaction equation: application to a pollution
  source problem,'' \emph{Inverse Problems}, vol.~21, no.~3, pp. 1121--1136,
  2005.

\bibitem{Komornik2005}
V.~Komornik and M.~Yamamoto, ``Estimation of point sources and applications to
  inverse problems,'' \emph{Inverse Problems}, vol.~21, pp. 2051--2070, 2005.

\bibitem{Albini2003451}
L.~Albini, P.~Burrascano, and S.~Fiori, ``A feasibility study for
  electromagnetic pollution monitoring by electromagnetic-source localization
  via neural independent component analysis,'' \emph{Neurocomputing}, vol.~55,
  no. 3–4, pp. 451 -- 468, 2003.

\bibitem{Chemical_plume}
S.~Pang and J.~Farrell, ``Chemical plume source localization,'' \emph{Part B:
  Cybernetics, IEEE Transactions on Systems, Man, and Cybernetics}, vol.~36,
  no.~5, pp. 1068--1080, 2006.

\bibitem{Following_RF}
A.~Wadhwa, U.~Madhow, J.~Hespanha, and B.~Sadler, ``Following an rf trail to
  its source,'' in \emph{49th Annual Allerton Conference on Communication,
  Control, and Computing (Allerton)}, 2011, pp. 580--587.

\bibitem{IslerTrackingFish}
P.~Tokekar, E.~Branson, J.~Vander~Hook, and V.~Isler, ``Tracking aquatic
  invaders: Autonomous robots for monitoring invasive fish,'' \emph{Robotics
  Automation Magazine, IEEE}, vol.~20, no.~3, pp. 33--41, 2013.

\bibitem{PSO}
J.~Kennedy and R.~Eberhart, ``Particle swarm optimization,'' in \emph{IEEE
  International Conference on Neural Networks}, vol.~4, 1995, pp. 1942--1948
  vol.4.

\bibitem{ModifiedPSO}
Y.~Shi and R.~Eberhart, ``A modified particle swarm optimizer,'' in \emph{IEEE
  International Conference on Evolutionary Computation Proceedings, IEEE World
  Congress on Computational Intelligence}, 1998, pp. 69--73.

\bibitem{Clerc_Constriction}
M.~Clerc and J.~Kennedy, ``The particle swarm - explosion, stability, and
  convergence in a multidimensional complex space,'' \emph{IEEE Transactions on
  Evolutionary Computation}, vol.~6, no.~1, pp. 58--73, 2002.

\bibitem{NeighborPSO}
K.~Veeramachaneni, T.~Peram, C.~Mohan, and L.~Osadciw,
  ``\BIBforeignlanguage{English}{Optimization using particle swarms with near
  neighbor interactions},'' in \emph{\BIBforeignlanguage{English}{Genetic and
  Evolutionary Computation}}.\hskip 1em plus 0.5em minus 0.4em\relax Springer
  Berlin Heidelberg, 2003, vol. 2723, pp. 110--121.

\bibitem{QuantumPSO}
M.~Pant, R.~Thangaraj, and A.~Abraham, ``A new quantum behaved particle swarm
  optimization,'' in \emph{Proceedings of the 10th Annual Conference on Genetic
  and Evolutionary Computation}.\hskip 1em plus 0.5em minus 0.4em\relax New
  York, NY, USA: ACM, 2008, pp. 87--94.

\bibitem{Vijay_digital}
V.~Kalivarapu, J.-L. Foo, and E.~Winer,
  ``\BIBforeignlanguage{English}{Improving solution characteristics of particle
  swarm optimization using digital pheromones},''
  \emph{\BIBforeignlanguage{English}{Structural and Multidisciplinary
  Optimization}}, vol.~37, no.~4, pp. 415--427, 2009.

\bibitem{PSOPugh2006}
J.~Pugh, L.~Segapelli, and A.~Martinoli,
  ``\BIBforeignlanguage{English}{Applying aspects of multi-robot search to
  particle swarm optimization},'' in \emph{\BIBforeignlanguage{English}{Ant
  Colony Optimization and Swarm Intelligence}}, ser. Lecture Notes in Computer
  Science.\hskip 1em plus 0.5em minus 0.4em\relax Springer Berlin Heidelberg,
  2006, vol. 4150, pp. 506--507.

\bibitem{PSOPugh2007}
J.~Pugh and A.~Martinoli, ``Inspiring and modeling multi-robot search with
  particle swarm optimization,'' in \emph{IEEE Swarm Intelligence Symposium},
  April 2007, pp. 332--339.

\bibitem{PSOJatmiko}
W.~Jatmiko, K.~Sekiyama, and T.~Fukuda, ``A pso-based mobile robot for odor
  source localization in dynamic advection-diffusion with obstacles
  environment: theory, simulation and measurement,'' \emph{Computational
  Intelligence Magazine, IEEE}, vol.~2, no.~2, pp. 37--51, May 2007.

\bibitem{PSOMarques}
L.~Marques, U.~Nunes, and A.~de~Almeida,
  ``\BIBforeignlanguage{English}{Particle swarm-based olfactory guided
  search},'' \emph{\BIBforeignlanguage{English}{Autonomous Robots}}, vol.~20,
  no.~3, pp. 277--287, 2006.

\bibitem{PSODerr}
K.~Derr and M.~Manic, ``Multi-robot, multi-target particle swarm optimization
  search in noisy wireless environments,'' in \emph{Human System Interactions,
  2009. HSI '09. 2nd Conference on}, May 2009, pp. 81--86.

\bibitem{PSODoctor}
S.~Doctor, G.~Venayagamoorthy, and V.~Gudise, ``Optimal pso for collective
  robotic search applications,'' in \emph{Evolutionary Computation, 2004.
  CEC2004. Congress on}, vol.~2, June 2004, pp. 1390--1395 Vol.2.

\bibitem{HerefordPSO2007}
J.~Hereford, M.~Siebold, and S.~Nichols, ``Using the particle swarm
  optimization algorithm for robotic search applications,'' in \emph{IEEE Swarm
  Intelligence Symposium}, April 2007, pp. 53--59.

\bibitem{PSO_Rui_AIM}
R.~Zou, V.~Kalivarapu, J.~Oliver, and S.~Bhattacharya, ``Swarm optimization
  techniques for multi-agent source localization,'' in \emph{IEEE/ASME
  International Conference on Advanced Intelligent Mechatronics (AIM)}, July
  2014, pp. 402--407.

\bibitem{PSO_Rui_MSC}
R.~Zou, M.~Zhang, V.~Kalivarapu, E.~Winer, and S.~Bhattacharya, ``Particle
  swarm optimization for source localization in environment with obstacles,''
  in \emph{IEEE Multi-Conference on Systems and Control}, October 2014, {T}o
  Appear.

\bibitem{rawnet11}
A.~Khanafer, S.~Bhattacharya, and T.~Ba\c{s}ar, ``Adaptive resource allocation
  in jamming teams using game theory,'' in \emph{7th International workshop on
  Resource Allocation and Cooperation in Wireless Networks(RAWNET)}, 2011.

\bibitem{PASPWEB2010}
J.~O. Smith, \emph{Physical Audio Signal Processing}.\hskip 1em plus 0.5em
  minus 0.4em\relax http://ccrma.stanford.edu/~jos/pasp, 2010, online book.

\bibitem{PSOPara}
Y.~Shi and R.~Eberhart, ``Parameter selection in particle swarm optimization,''
  in \emph{Evolutionary Programming VII}, ser. Lecture Notes in Computer
  Science.\hskip 1em plus 0.5em minus 0.4em\relax Springer Berlin Heidelberg,
  1998, vol. 1447, pp. 591--600.

\bibitem{kennedy2001swarm}
J.~Kennedy, R.~Eberhart, and Y.~Shi, \emph{Swarm Intelligence}, ser.
  Evolutionary Computation Series.\hskip 1em plus 0.5em minus 0.4em\relax
  Morgan Kaufmann Publishers, 2001.

\bibitem{SPSO}
D.~Bratton and J.~Kennedy, ``Defining a standard for particle swarm
  optimization,'' in \emph{IEEE Swarm Intelligence Symposium}, April 2007, pp.
  120--127.

\bibitem{SPSO2011}
M.~Zambrano-Bigiarini, M.~Clerc, and R.~Rojas, ``Standard particle swarm
  optimisation 2011 at cec-2013: A baseline for future pso improvements,'' in
  \emph{IEEE Congress on Evolutionary Computation (CEC)}, June 2013, pp.
  2337--2344.

\bibitem{clerc2010particle}
M.~Clerc, \emph{Particle Swarm Optimization}, ser. ISTE.\hskip 1em plus 0.5em
  minus 0.4em\relax Wiley, 2010.

\bibitem{bug_algorithm}
V.~J. Lumelsky and A.~Stepanov, ``Dynamic path planning for a mobile automaton
  with limited information on the environment,'' \emph{IEEE Transactions on
  Automatic Control}, vol.~31, no.~11, pp. 1058--1063, Nov 1986.

\bibitem{choset2005principles}
H.~Choset, \emph{Principles of Robot Motion: Theory, Algorithms, and
  Implementation}, ser. A Bradford book.\hskip 1em plus 0.5em minus 0.4em\relax
  Prentice Hall of India, 2005.

\bibitem{Dijkstra}
E.~Dijkstra, ``\BIBforeignlanguage{English}{A note on two problems in connexion
  with graphs},'' \emph{\BIBforeignlanguage{English}{Numerische Mathematik}},
  vol.~1, no.~1, pp. 269--271, 1959.

\end{thebibliography}
\end{document}